\documentclass[12pt]{article}
\usepackage{amssymb, amsfonts, amsmath, amsthm, url, graphicx}

\def\3{\subset }
\def\4{\subseteq }
\def\<{\left<}
\def\>{\right>}
\def\vsp{\vspace*{1,5mm}\\ }

\def\bit{\begin{itemize}}
\def\eit{\end{itemize}}
\def\3{\subset }
\def\4{\subseteq }

\def\0{\leqno}

\def\barr{\begin{array}}
\def\earr{\end{array}}
\def\dd{\displaystyle}

\def\Z{{\rlap{$\kern2pt{\rm Z}$}{\rm Z}\,}}

\def\frax{\dd\frac}

\title{\bf A note on subgroup commutativity degrees of finite groups}
\author{Marius T\u arn\u auceanu}
\date{2016/2017}

\begin{document}

\maketitle

\begin{abstract}
In this note we give some new results concerning the subgroup commutativity degree of a finite group $G$. These are obtained
by considering the minimum of subgroup commutativity degrees of all sections of $G$.
\end{abstract}



\section{Introduction}

In the last years there has been a growing interest in the use of
probability in finite group theory. One of the most important aspects
which have been studied is the probability that two elements of a
finite group $G$ commute. It is called the {\it commutativity degree}
of $G$ and is denoted by $d(G)$. Inspired by this concept, in \cite{8}
(see also \cite{9}) we introduced a similar notion for the subgroups
of $G$, called the {\it subgroup commutativity degree} (or the {\it subgroup
permutability degree}) of $G$. This quantity is defined by
$$\barr{lcl} sd(G)&=&\frax1{|L(G)|^2}\,\left|\{(H,K)\in L(G)^2\mid
HK=KH\}\right|=\vsp &=&\frax1{|L(G)|^2}\,\left|\{(H,K)\in
L(G)^2\mid HK\in L(G)\}\right|\earr$$ (where $L(G)$ denotes
the subgroup lattice of $G$) and it measures the probability that
two subgroups of $G$ commute, or equivalently the probability that
the product of two subgroups of $G$ be a subgroup of $G$.

We recall that for a finite group $G$ we have $sd(G)=1$ if and only if $G$
is an Iwasawa group, i.e. a nilpotent modular group (see [6, Exercise 3, p. 87]).
A complete description of these groups is given by a well-known Iwasawa's result
(see Theorem 2.4.13 of \cite{6}). In particular, we infer that $sd(G)=1$ for
all Dedekind groups $G$.

A well-known result by Gustafson \cite{3} concerning the commutativity degree
states that if $d(G)>5/8$ then $G$ is abelian, and we have $d(G)=5/8$ if
and only if $G/Z(G)\cong\mathbb{Z}_2\times\mathbb{Z}_2$. Note that the
similar problem for the subgroup commutativity degree does not have a solution,
i.e. there is no constant $c\in(0,1)$ such that if $sd(G)>c$ then $G$ is Iwasawa,
as shows Theorem 2.15 of \cite{1}.

In the following we will study this problem by replacing the condition "$sd(G)>c$" with
the stronger condition "$sd^*(G)>c$", where $$sd^*(G)={\rm min}\{sd(S) \mid S \mbox{ section of } G\}.$$It
was suggested by the fact that a $p$-group is modular if and only if each of its sections
of order $p^3$ is. Moreover, if a $p$-group is not modular then it contains a section
isomorphic to $D_8$, the dihedral group of order $8$, or to $E(p^3)$, the non-abelian group of order $p^3$
and exponent $p$ for $p>2$ (see Lemma 2.3.3 of \cite{6}). Note that a similar condition for the cyclic subgroup
commutativity degree led in \cite{11} to a criterion for a finite group to be an Iwasawa group. We will prove that
the condition $sd^*(G)>23/25$ implies the modularity for finite nilpotent groups $G$, and also that it implies
the solvability for arbitrary finite groups $G$. Then we will show the non-existence of a constant $c\in(0,1)$
such that if $sd^*(G)>c$ then $G$ is Iwasawa, extending the above mentioned result of Aivazidis.

Most of our notation is standard and will usually not be repeated here. Elementary notions
and results on groups can be found in \cite{7}. For subgroup lattice concepts we refer the
reader to \cite{6}.

\section{Main results}

\bigskip\noindent{\bf Theorem 1.} {\it Let $G$ be a finite nilpotent group such that $sd^*(G)>23/25$. Then $G$ is modular, and consequently an Iwasawa group.}
\bigskip

\noindent{\bf Proof.} Being nilpotent, $G$ can be written as a direct product of its Sylow subgroups $G_i$, $i=1,2,...,k$. Clearly, for each $i$ we have
$$sd^*(G_i)\geq sd^*(G)>\frax{23}{25}\,.$$Assume that $G_i$ is not modular. Then there is a section $S$ of $G_i$ such that $S\cong D_8$ or $S\cong E(p^3)$ for $p>2$. We can easily check that
$$sd(E(p^3))=\frax{3p^3+12p^2+16p+16}{(p^2+2p+4)^2}<\frax{23}{25}=sd(D_8)$$and therefore $sd(S)\leq 23/25$, a contradiction. Thus $G_i$ is modular, for all $i=1,2,...,k$, which implies that $G$ is itself a modular group.
\hfill\rule{1,5mm}{1,5mm}

\bigskip\noindent{\bf Remark.} The constant $23/25$ in Theorem 1 can be decreased for $p$-groups with $p>2$ by observing that such a group cannot have sections isomorphic to $D_8$. Thus, a finite $p$-group $G$ of odd order which satisfies $$sd^*(G)>\frax{3p^3+12p^2+16p+16}{(p^2+2p+4)^2}\0(1)$$\,is always an Iwasawa group. We also observe that $$\frax{3p^3+12p^2+16p+16}{(p^2+2p+4)^2}<\frax{3}{p}\,, \mbox{ for all } p$$\,and therefore the above statement remains true by replacing the condition (1) with the more elegant condition $$sd^*(G)\geq \frax{3}{p}\,.$$

In what follows we will study what can be said about an \textit{arbitrary} finite group $G$ satisfying $sd^*(G)>23/25$. A first answer is given by the following theorem.

\bigskip\noindent{\bf Theorem 2.} {\it Let $G$ be a finite group such that $sd^*(G)>23/25$. Then $G$ is solvable.}
\bigskip

\noindent{\bf Proof.} Assume that $G$ is not solvable. Then it contains a section isomorphic to one of the following groups:
\begin{itemize}
\item[-] ${\rm PSL}(2,p)$, where $p>3$ is a prime such that $5\mid p^2+1$;
\item[-] ${\rm PSL}(2,3^p)$, where $p\geq 3$ is a prime;
\item[-] ${\rm PSL}(2,2^p)$, where $p$ is a prime;
\item[-] ${\rm Sz}(2^p)$, where $p\geq 3$ is a prime;
\item[-] ${\rm PSL}(3,3)$.
\end{itemize}It is well-known that ${\rm PSL}(2,q)$ has a subgroup isomorphic to $D_{q+1}$ for $q$ odd and a subgroup isomorphic to $D_{2(q+1)}$ for $q=2^p$ (see e.g. \cite{2}). Then the first three groups above have a section isomorphic to $D_8$ or to $D_{2r}$ with $r\geq 3$ a prime. But $$sd(D_{2r})=\frax{7r+9}{(r+3)^2}<\frax{23}{25}\,,$$a contradiction. A similar contradiction is obtained for ${\rm Sz}(2^p)$ since it contains a subgroup of type $D_{2(2^p-1)}$. Finally, ${\rm PSL}(3,3)$ has a subgroup isomorphic to $A_4$ and $$sd(A_4)=\frax{16}{25}<\frax{23}{25}\,,$$contradicting again our hypothesis. This completes the proof.
\hfill\rule{1,5mm}{1,5mm}\bigskip

Next we try to see whether the condition $sd^*(G)>23/25$ (or the more general condition $sd^*(G)>c$) implies that $G$ is nilpotent. We start by providing an example of a Schmidt group $S_1$ of order $p^rq$ for which we are able to compute explicitly $sd^*(S_1)$. It also has the property that if $p$ and $q$ are suitably chosen, then  $sd^*(S_1)$ tends to $1$ when $p$ tends to infinity. This will be the main ingredient of the proof of Theorem 3.

\bigskip\noindent{\bf Example.} Let $S$ be a Schmidt group, i.e. a finite non-nilpotent group all of whose proper subgroups are nilpotent. By \cite{5} (see also \cite{4}) it follows that $S$ is a solvable group of order $p^mq^n$ (where $p$ and $q$ are different primes) with a unique Sylow $p$-subgroup $P$ and a cyclic Sylow $q$-subgroup $Q$, and hence $S$ is a semidirect product of $P$ by $Q$. Moreover, we have:
\begin{itemize}
\item[-] if $Q=\langle y\rangle$ then $y^q\in Z(S)$;
\item[-] $Z(S)=\Phi(S)=\Phi(P)\times\langle y^q\rangle$, $S'=P$, $P'=(S')'=\Phi(P)$;
\item[-] $|P/P'|=p^r$, where $r$ is the order of $p$ modulo $q$;
\item[-] if $P$ is abelian, then $P$ is an elementary abelian $p$-group of order $p^r$ and $P$ is a minimal normal subgroup of $S$;
\item[-] if $P$ is non-abelian, then $Z(P)=P'=\Phi(P)$ and $|P/Z(P)|=p^r$.
\end{itemize}We infer that $S_1=S/Z(S)$ is also a Schmidt group of order $p^rq$ which can be written as semidirect product of an elementary abelian $p$-group $P_1$ of order $p^r$ by a cyclic group $Q_1$ of order $q$ (note that $S_3$ and $A_4$ are examples of such groups). It is easy to see that $S_1$ does not contain subgroups of order $p^iq$ for $i=1,2,...,r-1$. Then $$L(S_1)=L(P_1)\cup\{Q_1^x \mid x\in S_1\}\cup\{S_1\}$$and so $$|L(S_1)|=a_{r,p}+p^r+1\,,$$where $a_{r,p}$ denotes the total number of subgroups of $P_1$. By \cite{10} the numbers $a_{r,p}$ can be written as $$a_{r,p}=f_r(p), \mbox{ where } f_r\in{\rm\mathbb{Z}[X]} \mbox{ and } {\rm deg}(f_r)=\left[r^2/4\right].$$Since $S_1/1$ is the unique non-abelian section of $S_1$, one obtains
$$sd^*(S_1)=sd(S_1).$$Let $Q_1^1,Q_1^2,...,Q_1^{p^r}$ be the conjugates of $Q_1$. Then the pairs of commuting subgroups of $S_1$ are:
\begin{itemize}
\item[-] $(X,Y)$ with $X,Y\leq P_1$,
\item[-] $(X,S_1)$ and $(S_1,X)$ with $X\leq P_1$,
\item[-] $(Q_1^i,1)$ and $(1,Q_1^i)$, $i=1,2,...,p^r$,
\item[-] $(Q_1^i,P_1)$ and $(P_1,Q_1^i)$, $i=1,2,...,p^r$,
\item[-] $(Q_1^i,S_1)$ and $(S_1,Q_1^i)$, $i=1,2,...,p^r$,
\item[-] $(Q_1^i,Q_1^i)$, $i=1,2,...,p^r$,
\item[-] $(S_1,S_1)$.
\end{itemize}It follows that\newpage
$$\barr{lcl} sd^*(S_1)&=&\frax{a_{r,p}^2+2a_{r,p}+7p^r+1}{(a_{r,p}+p^r+1)^2}=\vsp &=&\frax{1+\frax{2}{a_{r,p}}+\frax{7p^r}{a_{r,p}^2}+\frax{1}{a_{r,p}^2}}{1+\frax{p^{2r}}{a_{r,p}^2}+\frax{1}{a_{r,p}^2}+
\frax{2p^r}{a_{r,p}}+\frax{2}{a_{r,p}}+\frax{2p^r}{a_{r,p}^2}}\,.\earr\0(2)$$

\bigskip\noindent{\bf Theorem 3.} {\it There is no constant $c\in(0,1)$ such that if $sd^*(G)>c$ then $G$ is Iwasawa.}
\bigskip

\noindent{\bf Proof.} Let $(p_n)_{n\geq 1}$ and $(q_n)_{n\geq 1}$ be two strictly increasing sequences of primes such that the order $r_n$ of $p_n$ modulo $q_n$ is greater than $4$. It follows that $$\left[r_n^2/4\right]>r_n, \,\forall\, n\geq 1.\0(3)$$For every $n\geq 1$, let $G_n$ be a semidirect product of an elementary abelian $p_n$-group $P_n$ of order $p_n^{r_n}$ by a cyclic group of order $q_n$ generated by an element $x_n$ which permutes the elements of a basis of $P_n$ cyclically. Then $(G_n)_{n\geq 1}$ are Schmidt groups of order $p_n^{r_n}q_n$ such as $S_1$ in our example, and so (2) and (3) lead to $$\dd\lim_{n\to\infty}sd^*(G_n)=1,$$completing the proof.\hfill\rule{1,5mm}{1,5mm}
\bigskip

Inspired by the above results, we came up with the following conjecture.

\bigskip\noindent{\bf Conjecture 4.} {\it Let $G$ be a finite group such that $sd^*(G)>23/25$. Then $G$ is either an Iwasawa group or a Schmidt group.}

\bigskip\noindent{\bf Remark.} The above example also leads to two new classes of finite groups whose subgroup commutativity degree vanishes asymptotically. For $i=1,2$, let $(p^i_n)_{n\geq 1}$ and $(q^i_n)_{n\geq 1}$ be two strictly increasing sequences of primes such that the order $r^i_n$ of $p^i_n$ modulo $q^i_n$ is $2$ and $3$, respectively. Let $(G^i_n)_{n\geq 1}$ be the Schmidt groups of order $(p^i_n)^{r^i_n}q^i_n$ constructed as in the proof of Theorem 3. Then (2) implies that $$\dd\lim_{n\to\infty}sd^*(G^i_n)=\dd\lim_{n\to\infty}sd(G^i_n)=0,\, i=1,2,$$as desired.\bigskip
\newpage

Finally, we formulate another natural problem concerning our study.

\bigskip\noindent{\bf Open problem.} {\it Describe the structure of finite groups $G$ satisfying $sd^*(G)=23/25$.}
\bigskip

\bigskip\noindent {\bf Acknowledgements.} The author wish to thank Stefanos Aivazidis for fruitful discussions on this subject. He is also grateful to the reviewers for their useful suggestions to improve this paper.

\vspace*{5ex}\small

\hfill
\begin{minipage}[t]{5cm}
Marius T\u arn\u auceanu \\
Faculty of  Mathematics \\
``Al.I. Cuza'' University \\
Ia\c si, Romania \\
e-mail: {\tt tarnauc@uaic.ro}
\end{minipage}

\end{document}